\documentclass[10pt]{amsart}

\usepackage{amsmath}
\usepackage{amsthm}
\usepackage{amsopn}
\usepackage{amssymb}
\usepackage[all]{xy}
\usepackage{lscape,xcolor}
\usepackage{graphicx}
\usepackage{hyperref}
\usepackage{mathtools}
\usepackage{stmaryrd}
\usepackage{enumitem}

\newcommand{\eitem}{\stepcounter{equation}\item}

\allowdisplaybreaks[1]

\newcommand{\mc}[1]{\mathcal{#1}}
\newcommand{\ul}[1]{\underline{#1}}

\newcommand{\mr}[1]{\mathrm{#1}}

\newcommand{\mit}[1]{\mathit{#1}}

\newcommand{\bra}[1]{\langle #1 \rangle}

\newcommand{\td}[1]{\widetilde{#1}}

\newcommand{\ZZ}{\mathbb{Z}}

\newcommand{\QQ}{\mathbb{Q}}

\newcommand{\FF}{\mathbb{F}}
\newcommand{\GG}{\mathbb{G}}

\def \HF2{\mr{H}\FF_2}

\newcommand{\Sp}{\mr{Sp}}
\newcommand{\sS}{\mr{sSet}}
\newcommand{\sGp}{\mr{sGp}}

\newcommand{\sHopf}{\mr{sHopf}}
\newcommand{\sLie}{\mr{sLie}}

\newcommand{\Mod}{\mr{Mod}}
\newcommand{\RMod}{\mr{Mod}^{\mit{rt}}}
\newcommand{\LMod}{\mr{Mod}^{\mit{lt}}}
\newcommand{\Prim}{\mr{Prim}}

 \newtheorem{thm}[equation]{Theorem}

 \newtheorem*{thm*}{Theorem}
 \newtheorem*{cor*}{Corollary}
 \newtheorem*{lem*}{Lemma}
 \newtheorem*{prop*}{Proposition}

 \theoremstyle{definition}

 \newtheorem{rmk}[equation]{Remark}

\newtheorem*{defn*}{Definition}
\newtheorem*{ex*}{Example}
\newtheorem*{exs*}{Examples}
\newtheorem*{rmk*}{Remark}
\newtheorem*{claim*}{Claim}

\numberwithin{equation}{section}
\numberwithin{figure}{section}
\DeclareMathOperator{\Ext}{Ext}

\DeclareMathOperator{\Ind}{Ind}
\DeclareMathOperator{\Seq}{Seq}

\DeclareMathOperator{\gr}{gr}

\newcommand{\E}[2]{\prescript{#1}{#2}{E}}

\newcommand{\cHom}{\mc{H}\mit{om}}
\newcommand{\cLie}{\mc{L}\mit{ie}}
\newcommand{\cComm}{\mc{C}\mit{omm}}
\newcommand{\Id}{\mr{Id}}

\title{Unstable homotopy groups and Lie algebras}
\author{Mark Behrens}
\author{Connor Malin}

\begin{document}

\begin{abstract}
We survey the role of Lie algebras in the study of unstable homotopy groups.
\end{abstract}

\maketitle
\tableofcontents

\section{Introduction}

Let $X$ be a pointed space.
The Whitehead product
$$ [-,-]: \pi_{i}(X) \otimes \pi_{j}(X) \to \pi_{i+j-1}(X) $$
gives $\pi_*(X)$ the structure of a graded shifted Lie algebra.
This structure is most easily conceptualized by its relationship to the Samelson product, which is given by the commutator on the loop space  
$$ \bra{-,-}:\pi_i(\Omega X) \otimes \pi_j(\Omega X) \to \pi_{i+j}(\Omega X). $$
Samelson showed that under the isomorphism $\pi_{i+1}(X) \cong \pi_i(\Omega X)$, the two products agree up to a sign \cite{Samelson}. 

This Lie algebra structure is fundamentally unstable in nature --- there is no vestige of it in the context of stable homotopy groups.  It captures the difference between unstable and stable homotopy groups in a manner made precise by Curtis's lower central series \cite{Curtis}, Rector's mod $p$ lower central series \cite{Rector} and its relationship to simplicial restricted Lie algebras, and Quillen's differential graded Lie algebra model of unstable rational homotopy theory \cite{Quillen}.    

In this paper we will review these now classical ideas, and their more recent development in the context of Goodwillie calculus \cite{Goodwillie}.  Specifically, we will discuss the results of Arone, Ching, Taggart and the second author on Koszul duality and its interaction with Goodwillie calculus \cite{Ching}, \cite{AroneChing}, \cite{Espic}, \cite{MalinTaggart}, Konovalov's work on simplicial restricted Lie algebras \cite{Konovalov}, and the generalization of Quillen's Lie algebra model of rational homotopy theory to the unstable $v_n$-periodic context of Heuts, Rezk, and the first author \cite{Heuts}, \cite{BehrensRezk}. 

The recurring theme will be the following:
$$
\text{Unstable homotopy groups} = 
\left( \begin{array}{c} 
\text{stable homotopy groups} \\
+ \\
\text{Lie algebra information}
\end{array} \right)
$$
\subsection*{Conventions}
In an attempt to make the statements in this paper model independent, to free us from needing to specify hypotheses to express these statements in the proper derived context, and to allow us to present results spanning many decades and conventions in a unified and unambiguous language, we elect to make our statements on the level of $\infty$-categories.  In this sense they are valid in any of the equivalent models of $\infty$-categories the reader prefers \cite{Bergner}.  However, we note the Barwick-Kan model of relative categories provides the most convenient means of expressing classical results in the language of $\infty$-categories \cite{BarwickKan}.

We will denote the following $\infty$-categories by
\begingroup
\allowdisplaybreaks
\begin{align*}
\sS & = \text{the relative category of simplicial sets and weak homotopy equivalences} \\
& \quad \text{(a.k.a. \emph{spaces})} \\
\sS_* & = \text{pointed spaces} \\
\sS^{\ge n} & = \text{$(n-1)$-connected spaces} \\
\Sp & = \text{spectra}
\end{align*}
\endgroup
We will use $(-)^\vee$ to denote the Spanier-Whitehead dual.
Given an $\infty$-category $\mc{C}$, and objects $X,Y \in \mc{C}$, we let
$$ \mc{C}(X,Y) $$
denote the associated space of maps, and  
$$ [X,Y] = [X,Y]_{\mc{C}} $$
denote the corresponding set of homotopy classes of maps.  If $\mc{C}$ is a stable $\infty$-category, we will let
$$ \ul{\mc{C}}(X,Y) $$
denote the mapping spectrum.  

$p$ will always denote a prime number.  For elements $x_i$ of a Lie algebra $L$, we will let 
$$ [x_1, \ldots, x_k ] = [x_1, [x_2, \cdots [x_{k-1},x_k]\cdots]] $$
denote the iterated Lie bracket.

\subsection*{Acknowledgments}
The authors would like to thank Nikolay Konovalov, for his generous explanations of his work and related topics. The authors are also grateful to the referees for their helpful suggestions to improve this manuscript.

\section{Algebras and modules over operads}\label{sec:operads}

\subsection*{Symmetric sequences}
Fix a presentably symmetric monoidal stable $\infty$-category $(\mc{C},\otimes, 1_{\mc{C}})$, and let 
$$ \Seq_\Sigma(\mc{C}) $$
denote the $\infty$-category of \emph{symmetric sequences} in $\mc{C}$, whose objects are sequences
$$ \{ \mc{X}_i \in \mc{C}^{B\Sigma_i} \}_{i \ge 0}. $$ 
We will identify $\mc{C}$ with the full subcategory of $\Seq_\Sigma(\mc{C})$ consisting of sequences concentrated in degree $0$:
$$ X := \{ X, 0, 0, \ldots \}. $$
The $\infty$-category $\Seq_\Sigma(\mc{C})$ has a monoidal structure $\circ$ (see \cite[Sec.~4.1.2]{Brantner}) given by
$$ (\mc{X} \circ \mc{Y})_i := \bigoplus_{i = i_1 + \cdots + i_k} \Ind^{\Sigma_i}_{\Sigma_{i_1, \ldots, i_k}} \mc{X}_k \otimes \mc{Y}_{i_1} \otimes \cdots \otimes \mc{Y}_{i_k}, $$
where $\Sigma_{i_1, \ldots, i_k} \to \Sigma_i$ is the group of symmetries of the partition\footnote{Here, a partition $i = i_1 + \cdots + i_k$ is allowed to have $i_j = 0$.} $i = i_1 + \cdots + i_k$ and $\Ind$ is the left adjoint of the restriction along this homomorphism.  The unit of this monoidal structure is $1_*$, given by
$$ 1_* := \{ 0, 1_\mc{C}, 0, 0, \ldots \} \in \Seq_\Sigma(\mc{C}). $$
The $\infty$-category $\Seq_\Sigma(\mc{C})$ also possesses a symmetric monoidal structure $\otimes$ given by
\begin{equation}\label{eq:tensor}
 (\mc{X} \otimes \mc{Y})_i := \bigoplus_{i = i_1+i_2} \Ind^{\Sigma_i}_{\Sigma_{i_1} \times \Sigma_{i_2}} \mc{X}_{i_1} \otimes \mc{Y}_{i_2}.
 \end{equation}

\subsection*{Operads}
An \emph{operad in $\mc{C}$} is a monoid in $(\Seq_\Sigma(\mc{C}),\circ)$.
\begin{enumerate}[label = (\theequation)]   
\eitem We shall say that an operad $\mc{O}$ is \emph{reduced} if $\mc{O}_0 = 0$ and $\mc{O}_1 = 1$.  
\eitem Given an operad $\mc{O}$ in $\mc{C}$, we will let $\RMod_{\mc{O}}$/$\LMod_{\mc{O}}$ denote right/left modules over $\mc{O}$.
\eitem The $\otimes$-product of right $\mc{O}$-modules is again an $\mc{O}$-module, so $\otimes$ endows $\RMod_{\mc{O}}$ with a symmetric monoidal structure (see \cite{Fresse}). 
\eitem If an operad $\mc{O}$ is reduced, then the canonical map
$$ \mc{O} \to 1_* $$
is a map of operads, and in particular $1_*$ is both a left and right $\mc{O}$-module. 
\eitem An object $X \in \mc{C}$ gives rise to a symmetric sequence $X^{\otimes} \in \Seq_\Sigma(\mc{C})$ with
$$ ({X}^{\otimes})_i := X^{\otimes i}. $$
\eitem
An $\mc{O}$-coalgebra structure on $X$ is encoded in structure maps
$$ X \otimes \mc{O}_i \rightarrow X^{\otimes i}. $$
These maps induce a right $\mc{O}$-module structure on $X^{\otimes}$.
\eitem A left $\mc{O}$-module structure on $X \in \mc{C}$ (regarded as a symmetric sequence concentrated in degree $0$) is an $\mc{O}$-algebra structure on $X$.  This is because
$$ (\mc{O} \circ X)_0 \simeq \bigoplus_{0 = \underbrace{\scriptstyle{0+\cdots+0}}_{k}} \Ind^{\Sigma_0}_{\Sigma_k} \mc{O}_k \otimes X^{\otimes k} \simeq \bigoplus_k \mc{O}_k \otimes_{h\Sigma_k} X^{\otimes k}. $$
\end{enumerate}

\subsection*{Coendomorphism operads}
For objects $X,Y \in \mc{C}$, define a symmetric sequence of spectra $\cHom_{\mc{C}}(X,Y)$ by
$$ \mc{H}\mit{om}_{\mc{C}}(X, Y)_i := \ul{\mc{C}}(X,Y^{\otimes i}). $$
\begin{enumerate}[label=(\theequation)]
\eitem The symmetric sequence $\cHom_{\mc{C}}(X,X)$ admits a canonical operad structure (sometimes referred to as the \emph{coendomorphism operad}) \cite[3.4.1]{Fresse}\footnote{In \cite{Fresse}, the spectrum $\mc{H}om_{\mc{C}}(X,Y)$ corresponds to what Fresse refers to as the \emph{endomorphism module} $\mit{End}_{Y,X}$ taken in the opposite category $\mc{C}^{\mr{op}}$.}.
\eitem An $\mc{O}$-coalgebra structure on $X$ is a map of operads $\mc{O} \to \cHom_{\mc{C}}(X,X)$.
\eitem The symmetric sequence of spectra $\cHom_{\mc{C}}(X,Y)$ is canonically a right $\cHom_{\mc{C}}(Y,Y)$-module and a left $\cHom_{\mc{C}}(X,X)$-module \cite[8.1.1,8.1.3]{Fresse}\footnotemark[1].  
\eitem The \emph{$n$-sphere operad} is defined to be the coendomorphism operad in spectra
$$ \mc{S}^n := \cHom_{\Sp}(S^n, S^n). $$
Using the fact that $\mc{C}$ is tensored over spectra, we can define the \emph{$n$th suspension} $\sigma^n \mc{O}$ of an operad $\mc{O}$ to be the operad given by the \emph{levelwise tensor product}
$$ (\sigma^n \mc{O})_i := \mc{S}^n_i \otimes \mc{O}_i. $$
If $A$ is an $\sigma^n \mc{O}$-algebra, then $\Sigma^{n}A$ is an $\mc{O}$-algebra \cite[7.2.2]{LodayVallette} \cite[Prop.~6.20]{MalinOnePoint}.
\end{enumerate}

\subsection*{Koszul duality}
Ching originally defined Koszul duality of operads/modules in spectra using bar constructions \cite{Ching}.  Recently, Espic \cite{Espic} introduced a more conceptual homotopy invariant construction, which he showed was equivalent to Ching's.

Given a reduced operad of spectra $\mc{O}$, its \emph{Koszul dual} is defined to be the coendomorphism operad
$$ K(\mc{O}) := \cHom_{\RMod_\mc{O}}(1_*, 1_*). $$
Given $\mc{M} \in \RMod_{\mc{O}}$, its \emph{Koszul dual} is defined to be the right $K(\mc{O})$-module 
$$ K_{\mc{O}}(\mc{M}) := \cHom_{\RMod_{\mc{O}}}(\mc{M}, 1_*). $$
There are equivalences \cite[Sec.~3.3]{Espic}, \cite[Prop.~4.2.4]{MalinTaggart}
\begin{align*}
K(\mc{O}) & \simeq B(1_*, \mc{O}, 1_*)^{\vee}, \\
K_{\mc{O}}(\mc{M}) & \simeq B(\mc{M}, \mc{O}, 1_*)^\vee
\end{align*}  
where $B(-,-,-)$ denotes the two-sided monoidal bar construction.  Given an $\mc{O}$-coalgebra $X$, we define the spectrum of \emph{primitives} by
$$ \Prim_{\mc{O}}(X) := \ul{\Mod}^{\mit{rt}}_{\mc{O}}(1_*, X^{\otimes}) $$
It follows from the definition of $K(\mc{O})$ that $\Prim_{\mc{O}}(X)$ naturally has the structure of a $K(\mc{O})$-algebra.

\subsection*{The Lie operad}
Let $\cComm$ be the reduced commutative operad in spectra, given by the symmetric sequence
$$ \{ 0, S, S, S, \ldots \}. $$
Define the spectral Lie operad to be the shift of the Koszul dual
$$ \cLie := \sigma K(\cComm). $$
It is shown in \cite[Ex.~9.50]{Ching} that there is an isomorphism of operads
$$ H\ZZ_*\cLie \cong \cLie^\ZZ, $$
where $\cLie^\ZZ$ denotes the Lie operad in abelian groups.  For a commutative ring $k$, algebras over 
$$ \cLie^k := k \otimes \cLie^\ZZ $$ 
in $\Mod_k$ are Lie algebras over $k$. 

\section{The Goodwillie spectral sequence}

\subsection*{The Goodwillie tower}
Goodwillie calculus \cite{Goodwillie}, \cite[Ch.6]{Lurie} associates to a reduced functor  
between presentable pointed $\infty$-categories a \emph{Taylor tower} of degree $n$ polynomial approximations
$$ F \to \cdots \to P_n(F) \to \cdots \to P_1(F). $$
In the context where $F$ is the identity functor
$$ \Id: \sS_* \to \sS_*, $$
the fibers take the form \cite[Cor.~2.3]{Johnson}, \cite[Cor.~8.8]{Ching}
$$ D_n(\Id)(X) \simeq \Omega^\infty \sigma^{-1}\cLie_n  \otimes_{h\Sigma_n} \Sigma^\infty X^{\otimes n} $$
and for $X$ connected and nilpotent, the map
$$ X \to P_\infty(\Id)(X) := \varprojlim_n P_n(\Id)(X) $$
is an equivalence \cite[Sec.~5]{AroneKankaanrinta}.
It follows that for a connected nilpotent space there is a \emph{Goodwillie spectral sequence}
\begin{equation}\label{eq:gss}
 \E{gss}{}_1^{t,*}(X) = \pi_t \cLie(\Sigma^{-1}\Sigma^\infty X) \Rightarrow \pi_{t+1} X, 
\end{equation}
where 
$$ \cLie(Y) \simeq \bigoplus_n \cLie_n \otimes_{h\Sigma_n} Y^{\otimes n} $$
is the free spectral Lie algebra on a spectrum $Y$.  Note that the $E_1$-term is a Lie algebra, and the GSS converges to a Lie algebra, \emph{but it has not been proven that this spectral sequence is a spectral sequence of Lie algebras}.

\subsection*{The homotopy and homology of free spectral Lie algebras}
We are led to compute the homotopy groups of $\cLie(Y)$ for $Y \in \Sp$.  
We shall do this for the $p$-completions for every prime $p$.  

The homotopy groups of any bounded below $p$-complete spectrum $Z$ can be studied using the \emph{mod $p$ Adams spectral sequence}
$$ \Ext^{s,t}_{\mc{A}^{op}}(\FF_p, (H\FF_p)_*Z) \Rightarrow \pi_{t-s}Z. $$
Here $\mc{A}$ is the mod $p$ Steenrod algebra, whose dual action gives $(H\FF_p)_*Z$ an $\mc{A}^{op}$-module structure.
Thus the input needed to study the homotopy groups of the $p$-completion of $\sigma^{-1}\!\cLie(Y)$ is the homology $(H\FF_p)_* \sigma^{-1}\cLie(Y)$.

We first consider the case of $p = 2$.  Suppose that $L$ is a $2$-complete $\cLie$-algebra. Since
\begin{equation}\label{eq:Lie_2}
\cLie_2 \simeq S^{\sigma-1}
\end{equation}
where $\sigma$ is the sign representation, the $\cLie$-algebra structure gives a map
$$ (H\FF_2)_* (\Sigma^{\sigma-1} L^{\otimes 2})_{h\Sigma_2} \to (H\FF_2)_* L. $$
In addition to endowing $(H\FF_2)_*L$ with the structure of a graded Lie algebra, it also gives rise to \emph{Lie-Dyer-Lashof operations} \cite[Sec.~1.5]{Behrens}
\begin{equation}\label{eq:LieDyerLashof}
 \bar{Q}^i: (H\FF_p)_t L \to (H\FF_p)_{t+i} L 
 \end{equation}
which satisfy the \emph{allowablity conditions} \cite[Thm.~6.3]{OAC}
\begin{itemize}
\item $\bar{Q}^i x = 0$ if $i < |x|$.
\item $\bar{Q}^i x = [x,x]$ if $i = |x|$.
\item $[x, \bar{Q}^i y] = 0$.
\end{itemize}
The algebra of all such operations $\bar{\mc{R}}$ is subject to Lie-Adem relations \cite[Sec. 1.4]{Behrens} which give rise to a basis of admissible monomials
$$ \bar{Q}^{i_1} \cdots \bar{Q}^{i_\ell} $$
with $i_m > 2i_{m+1}$.  

Antol\'in Camarena \cite{OAC} showed that $(H\FF_2)_*\cLie(Y)$ is the free allowable $\bar{\mc{R}}$-Lie algebra on $(H\FF_2)_*Y$.  Specifically, if $\{x_j\}$ is a basis of $(H\FF_2)_*Y$, then $(H\FF_2)_*\cLie(Y)$ has a basis
\begin{equation}\label{eq:Liebasis}
 \bar{Q}^{i_1} \cdots \bar{Q}^{i_\ell}[x_{j_1}, \ldots , x_{j_k}] 
 \end{equation}
where the brackets range over a basis of the free graded Lie algebra over $\FF_2$ on the generators $\{x_j\}$, $i_m > 2i_{m+1}$, and $i_\ell > |x_{j_1}| + \cdots +|x_{j_k}|$.

For $p$ odd, Kjaer \cite{KjaerLie} constructed the odd primary analog of the Lie-Dyer-Lashof operations (\ref{eq:LieDyerLashof}), and he showed that $(H\FF_p)_*\cLie(Y)$ admits a basis analogous to (\ref{eq:Liebasis}).  However he was unable to determine the odd primary Lie-Adem relations.  In the case of the prime 2, they were deduced in \cite{Behrens} from a classical computation of the transfer
$$ (H\FF_2)_*B\Sigma_4 \to (H\FF_2)_*B\Sigma_2 \wr \Sigma_2 $$
due to Kahn and Priddy \cite{Priddy}.
Surprisingly, the formula for the odd primary analog of this transfer was unknown.  One interesting corollary of the work of \cite{Konovalov} (which we will discuss in Section~\ref{sec:restricted}) is that Konovalov is able to determine these odd primary Lie-Adem relations.  

\section{The mod $p$ lower central series and restricted Lie algebras}\label{sec:restricted}

\subsection*{The Rector spectral sequence} Adapting the work of Curtis \cite{Curtis} to the $p$-primary setting, Rector \cite{Rector} studied the spectral sequence associated to the the mod $p$-lower central series
$$ \Gamma^p_s G = \bra{[g_1,\ldots, g_i]^{p^j} \: : \: ip^j \ge s} \le G $$
of a simplicial group $G$, whose associated graded is a \emph{simplicial graded restricted Lie algebra over $\FF_p$}.
A graded restricted Lie algebra $L_*$ (over a field of characteristic $p$) is a graded Lie algebra which possesses an additional operation 
$$ \xi: L_t \to L_{pt} $$ 
which satisfies certain axioms (see \cite[6.1]{MilnorMoore}).  

Consider the following diagram of $\infty$-categories and functors.\footnote{Here, as in many other places in this paper, this commutative diagram of $\infty$-categories is obtained by observing there is a corresponding commutative diagram of relative categories, where we are deriving functors when necessary using functorial left/right approximations arising from model category structures.}
\begin{equation}\label{eq:resliediag}
\xymatrix{
\sS_* \ar[d]_{\td{\FF}_p} \ar[r]^{\Omega} & 
\sGp \ar[r]^{\Gamma^p_\bullet} \ar[d]_{\FF_p} &
\sGp^{fil} \ar[r]^{\gr_\bullet} &
\sLie_{\FF_p}^{\mit{gr.res}} \ar[d]^{V} 
\\
\mr{sCoAlg}^{\cComm}_{\FF_p} \ar[r]_{C} & 
\sHopf_{\FF_p} \ar[r]_{I^\bullet} &
\sHopf^{fil}_{\FF_p} \ar[r]_{\gr_\bullet} & 
\sHopf^{\mit{gr.prim}}_{\FF_p}   
}
\end{equation}
where
\begin{align*}
\sGp & = \text{simplicial groups} \\
\sLie_{\FF_p}^{\mit{gr.res}} & = \text{simplicial graded restricted Lie algebras over $\FF_p$} \\
\mr{sCoAlg}^{\cComm}_{\FF_p} & = \text{(non-counital) simplicial cocommutative coalgebras over $\FF_p$} \\
\sHopf_{\FF_p} & = \text{simplicial cocommutative Hopf algebras over $\FF_p$} \\
\sHopf^{\mit{gr.prim}}_{\FF_p} & = \text{simplicial cocommutative primitively generated graded Hopf algebras}
\end{align*}
and
\begin{align*}
\Omega X & = \text{the Kan loop group of a simplicial set $X$} \\
\FF_p X & = \text{the free simplicial $\FF_p$-module of a simplicial set $X$} \\
\td{\FF}_p X & = \text{the free reduced simplicial $\FF_p$-module of a pointed simplicial set $X$} \\
V(L) & = \text{the universal enveloping algebra (of a graded restricted Lie algebra $L$)} \\
C(A) & = C(\FF_p, \FF_p \oplus A, \FF_p) \text{, the cobar construction on a non-counital coalgebra $A$ } \\
& \quad \text{(where $\FF_p \oplus A$ denotes the coaugmented counital coalgebra associated to $A$)} \\
I^\bullet A & = \text{the filtration of $A$ given by powers of the augmentation ideal} \\
\mc{C}^{fil} & = \text{filtered objects of $\mc{C}$} \\
\gr_\bullet & = \text{the associated graded of a filtered object}
\end{align*}
The left-hand square of (\ref{eq:resliediag}) commutes when restricted to $\sS_*^{\ge 2}$ by the convergence of the Eilenberg-Moore spectral sequence. The right-hand rectangle of (\ref{eq:resliediag}) is shown to commute in \cite[Lem.~10.1]{BousfieldCurtis}.  The universal enveloping algebra functor $V$ in (\ref{eq:resliediag}) is an equivalence by \cite[Thm.~6.11]{MilnorMoore}, where the inverse functor is given by taking primitives 
$$ \Prim: \sHopf_{\FF_p}^{\mit{gr.prim}} \to \sLie_{\FF_p}^{\mit{gr.res}}. $$
Because the Kan loop group is level-wise free, the image of $X$ under the various functors of (\ref{eq:resliediag}) is given by
\[ 
\xymatrix{
X \ar@{|->}[r] & 
\Omega X \ar@{|->}[r] \ar@{|->}[d] &
\Gamma^p_\bullet \Omega X \ar@{|->}[r] &
\cLie^{r}(\Sigma^{-1} \td{\FF}_p X) \ar@{|->}[d]  
\\
& \FF_p \Omega X \ar@{|->}[r] &
I^\bullet \FF_p \Omega X \ar@{|->}[r] &
T(\Sigma^{-1}\td{\FF}_p X) 
}
\]
where $\cLie^{r}$ denotes the free restricted Lie algebra and $T$ denotes the tensor algebra. Note that each of these carries a natural grading with $\Sigma^{-1}\td{\FF}_p X$ in degree 1.

The filtration $\Gamma^p_\bullet \Omega X$ gives rise to the Rector spectral sequence
$$ \E{rss}{}_1^{t,*}(X) = \pi_t \cLie^{r}(\Sigma^{-1}\td{\FF}_p X) \Rightarrow \pi_{t+1} X^\wedge_p $$
which converges for $X$ simply connected \cite[Thm.~4.1]{Rector}, \cite[p.109]{BousfieldKan}.  

\begin{rmk}
The spectral sequence associated to the filtration $I^\bullet \FF_p \Omega X$ is the Eilenberg-Moore spectral sequence.  Thus, by (\ref{eq:resliediag}), the Hurewicz homomorphism induces a map from the Rector spectral sequence to the Eilenberg-Moore spectral sequence.
\end{rmk}

\subsection*{The homotopy of free simplicial restricted Lie algebras} In order to compute the $E_1$-term of the Rector spectral sequence, we observe that the homotopy groups of a simplicial restricted Lie algebra $L$ over $\FF_p$ have algebraic structure \cite[Prop.~8.3]{BousfieldCurtis}.  The \emph{restricted} structure arises from a factorization of the Lie algebra structure maps through maps
$$ \cLie^{\FF_p}_n \otimes^{h\Sigma_n} L^{\otimes n} \to L. $$ 
This endows $\pi_*L$ with the structure of a graded restricted Lie algebra.  Furthermore we get $\lambda$-operations coming from the mod $p$ cohomology of $\Sigma_p$.  For simplicity, we restrict attention to the case where $p = 2$.  In this case, it follows from (\ref{eq:Lie_2}) that
the map
$$ \cLie^{\FF_2}_2 \otimes^{h\Sigma_2} L^{\otimes 2} \to L $$ 
induces operations
$$ \lambda_i: \pi_t L \to \pi_{t+i} L $$
for $i \ge 0$, which satisfy the \emph{instability conditions} 
\begin{itemize}
\item $x \lambda_i = 0$ if $i > |x|$. 
\item $x \lambda_i = \xi(x)$ if $i = |x|$.
\item $[x, y\lambda_i] = 0$ if $i < |y|$.
\end{itemize}
The algebra of all such operations $\Lambda$ is subject to Adem relations \cite{6A}, which give rise to a basis of admissible monomials
$$ \lambda_{i_1} \cdots \lambda_{i_\ell} $$
with $2i_m \ge i_{m+1}$.  The $\Lambda$-algebra is Koszul dual to the Steenrod algebra $\mc{A}$ \cite{Priddy} and as such possesses a differential $d$ such that
$$ H^*(\Lambda) = \Ext_{\mc{A}}(\FF_2, \FF_2). $$

Bousfield and Curtis \cite[Thm.~7.3]{BousfieldCurtis} showed that for a simplicial $\FF_2$-module $Y$, the homotopy groups $\pi_* \cLie^{r}(Y)$ form the free unstable $\Lambda$-Lie algebra on $\pi_* Y$.  Specifically, if $\{x_j\}$ is a basis of $\pi_* Y$, then $\pi_* \cLie^{r}(Y)$ has a basis
$$ [x_{j_1}, \ldots, x_{j_k}]\lambda_{i_1} \cdots \lambda_{i_\ell} $$ 
where the brackets range over a basis of the free graded Lie algebra over $\FF_2$ on the generators $\{x_j\}$, $2i_m \ge i_{m+1}$, and $i_1 < |x_{j_1}| + \cdots +|x_{j_k}|$.

The Rector spectral sequence is a spectral sequence of Lie algebras.  The first potentially non-trivial differential on Lie algebra generators is given by
the formula \cite[Thm.~12.1]{BousfieldCurtis}
\begin{equation}\label{eq:rssdiff}
 d^{rss}_{2^\ell}(\sigma^{-1} x \cdot \lambda_I) = \sum [\sigma^{-1} x_i',\sigma^{-1}x_i]\lambda_I + \sigma^{-1}x\cdot d\lambda_I + \sum_j \sigma^{-1}x\mr{Sq}^j_* \cdot \lambda_{j-1}\lambda_I 
\end{equation}
for $x \in (\td{H\FF}_2)_*X$ with $\Delta(x) = \sum_{i}x_i' \otimes x_i''$.  It is shown in \cite[Thm.~5.2]{BousfieldCurtis} that for sufficiently nice spaces, the Rector spectral sequence is isomorphic to the unstable Adams spectral sequence after re-indexing.

\subsection*{The algebraic Goodwillie spectral sequence}

Konovalov \cite{Konovalov} related the Rector spectral sequence to the Goodwillie spectral sequence.  Specifically, he showed that the \emph{algebraic Goodwillie spectral sequence} associated to the Goodwillie tower of the functor\footnote{Technically, Konovalov studied the case of restricted Lie algebras over the algebraic closure $\bar{\FF}_p$.  This was so he could use the action of the units of $\bar{\FF}_p$ to prove degeneration results --- his results then carry over to $\FF_p$.} 
$$ \cLie^{r}(\Sigma^{-1} \td{\FF}_p (-)): \sS_* \to \sLie^{\mit{res}}_{\FF_p} $$
takes the form
\begin{equation}\label{eq:agss}
 \E{agss}{}_1 = (H\FF_p)_* \cLie(\Sigma^{-1}\Sigma^\infty X) \otimes \Lambda \Rightarrow \pi_* \cLie^r(\Sigma^{-1}\td{\FF}_p X). 
 \end{equation}
The spectral sequence (\ref{eq:agss}) is a spectral sequence of Lie algebras, and Konovalov \cite[Rmk.~8.3.7]{Konovalov} showed that the spectral sequence is \emph{completely determined} by explicit differentials on Lie algebra generators given by formulas discovered by Lin \cite{Lin} in his proof of the algebraic Kahn-Priddy theorem.

The conjecture is that the algebraic Goodwillie spectral sequence fits into a ``commuting square'' of spectral sequences\footnote{By ``commuting square'' we mean that the square of spectral sequences arises from a bifiltered object.}
\[ 
\xymatrix{
(H\FF_p)_* \cLie(\Sigma^{-1}\Sigma^\infty X) \otimes \Lambda \ar@{=>}[r]^-{agss} \ar@{=>}[d]_{ass} & 
\pi_* \cLie^r(\Sigma^{-1}\td{\FF}_p X) \ar@{=>}[d]^{rss}
\\
\pi_* \cLie(\Sigma^{-1}\Sigma^\infty X)^\wedge_p \ar@{=>}[r]_{gss}  &
\pi_{*+1}X^\wedge_p 
} \]
Such a commuting square would allow for the lifting of AGSS differentials to GSS differentials in a manner similar to \cite{Miller}.

\section{Lie algebra models of rational homotopy theory}

\subsection*{Rational homotopy theory}

Quillen famously showed that simply connected rational homotopy theory can be modeled by simplicial Lie algebras over $\QQ$ \cite{Quillen}.  He accomplished this by observing that diagram (\ref{eq:resliediag}) simplifies to a diagram of equivalences of $\infty$-categories (the equivalence $(*)$ was also studied by Sullivan \cite{Sullivan}).
\begin{equation}\label{eq:Quillen}
\xymatrix{
(\sS_*)^{\ge 2}_\QQ \ar|{(*)}[d]^{\: \: \simeq}_{\td{\QQ} \: \:} \ar[r]_\simeq^{\Omega} & 
\sGp^{\ge 1}_\QQ  \ar[d]^\simeq_{\QQ}
\\
(\mr{sCoAlg}^{\cComm}_{\QQ})^{\ge 2} \ar[r]^\simeq_{C} \ar[d]_{\Prim_{\cComm}}^\simeq & 
\sHopf^{\mit{conn}}_{\QQ} \ar[d]|{\simeq}^{\: \Prim}
\\
\mr{Alg}_{\sigma^{-1}\cLie}(\Sp_\QQ)^{\ge 2} \ar[r]^-\simeq_-{\Sigma^{-1}} &
\sLie^{\ge 1}_\QQ \ar@/^/[u]^{U}
}
\end{equation}
Here, $(\sS_*)^{\ge 2}_\QQ$, $\sGp^{\ge 1}_\QQ$ refer to the $\infty$-categories associated to the relative categories of simply connected simplicial sets (respectively, connected simplicial groups) and rational equivalences, $U$ refers to the universal enveloping algebra, and $\Prim_{\cComm}$ is the derived primitives construction described in Section~\ref{sec:operads}. We shall let $(-)_\QQ$ denote the associated localization functor 
$$ (-)_\QQ : \sS_*^{\ge 2} \to (\sS_*)^{\ge 2}_\QQ. $$

\subsection*{The rational Goodwillie tower}

In the rational case, for $X$ simply connected, the Goodwillie spectral sequence is the spectral sequence obtained from the bracket-length filtration on $\Prim \QQ \Omega X$, and takes the form \cite[Sec.~8.1]{Walter}
$$ \cLie^\QQ(\Sigma^{-1}(\td{H\QQ})_*X) \Rightarrow \pi_{*+1}X_\QQ. $$
In this case the spectral sequence is known to be a spectral sequence of Lie algebras.  The $d_1$-differential is determined by its effect on Lie algebra generators:
for $x \in (\td{H\QQ})_*X$ with $\Delta(x) = \sum x'_i \otimes x''_i$ this differential is given by
\cite[Apx~B]{Quillen}
$$ d^{gss}_1(\sigma^{-1}x) = \frac{1}{2}\sum_i (-1)^{|x_i'|}[\sigma^{-1} x_i',\sigma^{-1} x_i'']. $$
Thus if $X$ is of finite type, the $E_1$-page is the Harrison complex associated to the ring $(H\QQ)^*X$, and the $E_2$-page is its Andre-Quillen cohomology.

\section{Lie algebra models of unstable $v_n$-periodic homotopy theory}

\subsection*{The Bousfield-Kuhn functor}
Let $K(n)$ denote height $n$ Morava $K$-theory (see, for example, \cite{RavenelOrange}).  Recall that a $p$-local finite complex $X$ is called \emph{type $n$} if it is $K(n-1)$-acyclic, and $K(n) \otimes X \not\simeq 0$.
The periodicity theorem of Hopkins-Smith \cite{HopkinsSmith} implies that a $p$-local finite complex $V$ of type $n$ admits an asymptotically unique \emph{$v_n$-self-map}: a $K(n)$-equivalence 
$$ v: \Sigma^{t+N} V \to \Sigma^t V $$
for $t \gg 0$.
The \emph{unstable $v_n$-periodic homotopy groups (with coefficients in $V$)} of a pointed space $X$ are defined to be\footnote{The operator $v$ only acts on $[\Sigma^* V, X]$ for $* \ge t$.  We then can extend $v_n^{-1}\pi_*(X;V)$ to be  $\ZZ$-graded using $v_n$-periodicity.} 
$$ v_n^{-1}\pi_*(X;V) :=  v^{-1}[\Sigma^* V, X]. $$
The corresponding stable $v_n$-periodic homotopy groups  
$$ v_n^{-1}\pi_*^s(X; V) := \varprojlim_k v_n^{-1} \pi_{*+k}(\Sigma^k X; V) $$
are the homotopy groups of the telescope
$$ v^{-1}V^\vee \otimes \Sigma^{\infty} X. $$
Thus the stable $v_n$-periodic homotopy type of $X$ is encoded in the Bousfield localization 
$$ (\Sigma^\infty X)_{T(n)} \in \Sp_{T(n)} $$
where $T(n) := v^{-1}V^\vee$ (this localization is independent of the choice of $V$ and $v$-self map).

The \emph{Bousfield-Kuhn functor} \cite{KuhnBK}
$$ \Phi_n: \sS_* \to \Sp_{T(n)} $$ 
encodes these unstable $v_n$-periodic homotopy groups, in the sense that there are natural isomorphisms
$$ \pi_* \Phi_n(X) \otimes V^{\vee} \cong v_n^{-1}\pi_*(X;V). $$
The \emph{completed unstable $v_n$-periodic homotopy groups} are defined to be
$$ v_n^{-1} \pi^{\wedge}_*(X) := \pi_* \Phi_n(X). $$

\subsection*{A generalization of rational homotopy theory}
Let $v_n^{-1}\sS_*$ denote the $\infty$-category obtained by inverting the $v_n^{-1}\pi^{\wedge}_*$-isomorphisms.  Heuts \cite{Heuts} showed that $\Phi_n(X)$ canonically admits the structure of a $\sigma^{-1} \cLie$ algebra which is compatible with the Whitehead product on homotopy groups, and proved that the induced functor
$$ \Phi_n: v_n^{-1}\sS_* \xrightarrow{\simeq} \mr{Alg}_{\sigma^{-1}\cLie}(\Sp_{T(n)}) $$
is an equivalence of $\infty$-categories.

While this result gives a fantastic generalization of Quillen's simplicial Lie model of rational homotopy theory, it tells us nothing about the homotopy type of $\Phi_n(X)$.  To that end, one can imitate Quillen's approach to the rational case.
It is shown in \cite[Rmk.~5.13]{Heuts} (see also \cite[Sec.~6]{BehrensRezk}) that the diagram
\[
\xymatrix{v_n^{-1}\sS_* \ar[dd]^{\simeq}_{\Phi_n} \ar[dr]^{\Sigma^\infty_{T(n)}} \\
& \mr{CoAlg}_{\cComm}(\Sp_{T(n)}) \ar[dl]^{\: \: \Prim_{\cComm}} \\
\mr{Alg}_{\sigma^{-1}\cLie}(\Sp_{T(n)})
}
\] 
is lax commutative in the sense that there is a natural transformation called the \emph{comparison map}
\begin{equation}\label{eq:comparison}
 c_X: \Phi_n(X) \to \Prim_{\cComm} \Sigma^{\infty}_{T(n)} X. 
\end{equation} 

It is natural to consider the interaction of the Bousfield-Kuhn functor $\Phi_n$ with the Goodwillie tower.  
We shall say the a space $X$ is \emph{$\Phi_n$-good} if the map
\begin{equation}\label{eq:phigood}
 \Phi_n(X) \to P_\infty({\Phi}_n)(X)
 \end{equation}
is an equivalence.

One of the main results of \cite{Heuts} is 
\begin{thm}[Heuts]\label{thm:heuts}
The comparison map (\ref{eq:comparison}) is an equivalence for $X$ which are $\Phi_n$-good. 
\end{thm}
Theorem~\ref{thm:heuts} improved upon the main result of \cite{BehrensRezk}, which showed that if $X$ is finite with (\ref{eq:phigood}) a $K(n)$-equivalence (i.e. $X$ is \emph{$\Phi_{K(n)}$-good}), then the comparison map (\ref{eq:comparison}) is a $K(n)$-equivalence.  If $n = 1$, the validity of the telescope conjecture \cite[Prop.~4.2]{BousfieldLoc} implies $\Phi_1(X)$ is $K(1)$-local.  The telescope conjecture has been shown to be false for $n > 1$ \cite{BHLS}.

Arone and Ching discovered yet another approach to Theorem~\ref{thm:heuts} in the case where $X$ is finite, assuming certain results about Koszul duality of right modules, which was described in \cite[Sec.~9]{BehrensRezkSurvey}.  These Koszul duality results have now been proven \cite{MalinTaggart}, and in the next two subsections we will proceed to give a concise recapitulation of the Arone-Ching approach to Theorem~\ref{thm:heuts}.

\subsection*{Koszul duality and calculus}
Goodwillie \cite{Goodwillie} associated to a functor
$$ F : \sS_* \to \Sp $$
a symmetric sequence $\partial_*(F)$ of spectra, so that the fibers of its Taylor tower $\{P_k(F)\}$ are given by
$$ \partial_k(F) \otimes_{h\Sigma_k} \Sigma^\infty X^{\otimes k}. $$
Define the \emph{Koszul dual derivatives} to be the spectrum of natural transformations
$$ \partial^k(F) := \ul{\mr{Nat}}_X(F(X),\Sigma^\infty X^{\otimes k}). $$
The diagonal of $X$ induces a right $\cComm$-module structure on $\partial^*(F)$.
Using the Yoneda lemma, there is a natural transformation
$$ F(X) \to \ul{\Mod}^{\mit{rt}}_{\cComm}(\partial^*(F), \Sigma^\infty X^{\otimes}) $$
which gives an approximation of $F(X)$.

The derivatives $\partial_*(F)$ were shown in \cite[Thm.~4.2.8]{AroneChing} to possess a right $\sigma^{-1}\cLie$-module structure. 
The reason that we refer to $\partial^*(F)$ as the Koszul dual derivatives of $F$ is that if each $\partial_i F$ is a finite spectrum,
there is an equivalence \cite[Defn~4.2.7, Thm.~4.2.8]{AroneChing} of right $\cComm = K(\sigma^{-1}\cLie)$-modules 
\begin{equation}\label{eq:dualderivatives}
 \partial^*(F) \simeq K_{\sigma^{-1}\cLie}(\partial_*(F)).
\end{equation}   
It follows from the results of \cite{MalinTaggart} that if $X$ and all of the derivatives $\partial_k(F)$ are finite, then Koszul duality gives an equivalence
$$  \ul{\Mod}^{\mit{rt}}_{\cComm}(\partial^*(F), \Sigma^\infty X^{\otimes})  
\simeq \ul{\Mod}^{\mit{rt}}_{\sigma^{-1}\cLie}(\partial_*(\Sigma^\infty \sS_*(X,-)), \partial_*(F)) =: \Psi(F)(X)
$$
where 
$$
\Psi(F)(X) \simeq  \varprojlim_k \Psi_k(F) = \varprojlim_k \ul{\Mod}^{\mit{rt}}_{\sigma^{-1}\cLie}(\partial_{\le k}(\Sigma^\infty \sS_*(X,-)), \partial_{\le k}(F))
$$
is the \emph{fake Taylor tower} of \cite{AroneChingfake}.\footnote{In \cite{AroneChing},\cite{AroneChingfake}, the notation $\Phi_k(F)$ is used, but we instead use $\Psi_k$ to avoid conflict with the notation for the Bousfield-Kuhn functor.}

There is a map from the Taylor tower to the fake Taylor tower, giving a diagram of fiber sequences \cite[Rmk.~4.2.27]{AroneChing}
\begin{equation}\label{eq:fake}
\xymatrix{
\partial_k(F) \otimes_{h\Sigma_k} \Sigma^\infty X^{\otimes k} \ar[r] \ar[d]_{N} &  
P_k(F)(X) \ar[d] \ar[r] & P_{k-1}(F)(X) \ar[d] 
\\
\partial_k(F) \otimes^{h\Sigma_k} \Sigma^\infty X^{\otimes k} \ar[r] &  
\Psi_k(F)(X) \ar[r] & \Psi_{k-1}(F)(X) 
}
\end{equation}
where $N$ is the norm map.

\subsection*{Proof of Theorem~\ref{thm:heuts}}  We may now explain how the theory of the previous subsection specializes in the case of $F = {\Phi}_n$ to prove Theorem~\ref{thm:heuts} in the case where $X$ is finite.  The structure of the Goodwillie tower of $\Phi_n$ will be determined by the following theorem, which may also be deduced from \cite[Thm.~2.11]{Heuts}.  

\begin{thm}
    There is a canonical equivalence of right $\sigma^{-1}\cLie $ modules $$\partial_\ast \Phi_n \simeq (\sigma^{-1}\cLie)_{T(n)}.  $$
\end{thm}

\begin{proof}
    By \cite[Thm.~1.1]{KuhnBK}, there is a natural equivalence
    $$\Phi_n \Omega^\infty \simeq \mathrm{Id}_{\Sp_{T(n)}}$$
    By applying the generalized chain rule \cite[Thm.~A]{BlansBlom}, we can construct a chain of equivalences
    $$ (1_*)_{T(n)} \simeq \partial_\ast \Phi_n \Omega^\infty \simeq B(\partial_\ast \Phi_n,\sigma^{-1}\cLie,1_* )_{T(n)} $$
    Taking cobar constructions over $\cComm$ and applying \cite[Proposition 6.1]{ChingBarCobar}, yields the claimed equivalence at the level of symmetric sequences. Hence, $T(n)$-locally the right module $\partial_\ast( \Phi_n)$ satisfies the finiteness assumption needed to apply \cite[Theorem A]{MalinTaggart} and to deduce the above equivalence can be upgraded to an equivalence of right $\sigma^{-1}\cLie$-modules, proving the result.
\end{proof}

It follows from (\ref{eq:dualderivatives}) that the dual derivatives are given by 
$$ \partial^*({\Phi}_n) \simeq (1_*)_{T(n)} $$
Since $\sigma^{-1}\cLie_k$ is level-wise finite, we may apply Koszul duality for right modules \cite{MalinTaggart} to deduce that there is an equivalence
$$ \Psi({\Phi}_n)(X) \simeq \ul{\Mod}^{\mit{rt}}_{\cComm}(1_*, \Sigma_{T(n)}^{\infty}X^{\otimes}) = \Prim_{\cComm}(\Sigma^\infty_{T(n)} X). $$
Finally, since Kuhn proved norm maps in $\Sp_{T(n)}$ are equivalences \cite[Thm.~1.5]{Kuhn}\footnote{Kuhn proves that for any finite $G$,  $((S_{T(n)})^{tG})_{T(n)} \simeq 0$; since for any $T(n)$-local Borel $G$-spectrum $Y$, the Tate spectrum $(Y^{tG})_{T(n)}$ is a module over this spectrum, it must be trivial.}, it follows from (\ref{eq:fake}) that there is a natural equivalence
$$ P_\infty(\Phi_n)(X) \xrightarrow{\simeq} \Psi({\Phi}_n) (X). $$
We deduce that the comparison map (\ref{eq:comparison}) may be identified with the composite
$$ \Phi_n(X) \to P_\infty({\Phi}_n)(X) \xrightarrow{\simeq} \Prim_{\cComm}(\Sigma^\infty_{T(n)} X). $$
Theorem~\ref{thm:heuts} follows.

\subsection*{The $v_n$-periodic Goodwillie spectral sequence}

The Taylor tower for ${\Phi}_n$ gives rise to the $v_n$-periodic Goodwillie spectral sequence (which converges when $X$ is $\Phi_n$-good)
$$ v^{-1}_n\E{\mit{gss}}{}^{t,*}_1(X) = \pi_t \cLie(\Sigma^{-1} \Sigma^{\infty} X)_{T(n)} \Rightarrow v_n^{-1}\pi^{\wedge}_{t+1} (X). $$
Arone and Mahowald \cite[Thm.~0.1, Thm.~4.1]{AroneMahowald} showed that in the case where $X = S^d$ (and $d$ is odd if $p$ is odd), $v^{-1}_n\E{gss}{}^{t,k}(S^d) = 0$ unless $k = p^i \le p^n$, and they use this to prove that spheres are $\Phi_n$-good.
The $v_1$-periodic GSS was computed for $S^d$ by Mahowald \cite{MahowaldEHP} for $p = 2$ and Thompson \cite{Thompson} for $p$ odd.

As $T(n)$-local homotopy groups are largely incomputable at present for $n > 1$, one may alternatively consider the $K(n)$-local Goodwillie spectral sequence
$$ \E{\mit{gss}}{K(n)}^{t,*}_1 = \pi_t \cLie(\Sigma^{-1} \Sigma^{\infty} X)_{K(n)} \Rightarrow \pi_{t+1} \Phi_{n}(X)_{K(n)}. $$
The $K(2)$-local GSS for $S^3$ and $p \ge 5$ was computed by Wang in \cite{Wang}.

In general, 
the homotopy groups of the $K(n)$-localization of a spectrum $Z$ may be computed by its $K(n)$-local Adams-Novikov spectral sequence, which by the Morava change of rings theorem \cite{Morava}
takes the form
$$ H^s_c(\GG_n; (E_n)_t Z) \Rightarrow \pi_{t-s} Z_{K(n)}. $$
Here $(E_n)_*Z$ is the \emph{(completed) Morava $E$-homology}, and $\GG_n$ is the $n$th \emph{(extended) Morava stabilizer group}.
Thus the input needed to study the $K(n)$-local GSS is 
$(E_n)_*\cLie(Z)$.  

The Morava $E$-theory of $\cLie(Z)$ was computed by Brantner \cite{Brantner} in the case where $(E_n)_*Z$ is flat over $(E_n)_*$.  We briefly summarize his result.  Let $\Delta$ denote the Dyer-Lashof algebra for Morava $E$-theory, which acts on the Morava $E$-cohomology of any space.  The algebra $\Delta$ was shown by Rezk to be Koszul \cite{Rezk}.  Define the algebra of \emph{Hecke operations} $\mc{H}^{\cLie}$ to be the twisted Koszul dual algebra of $\Delta$ (a twisted version of Priddy's Koszul dual construction \cite{PriddyKoszul} defined in \cite[Defn.~4.3.1]{Brantner}).  For simplicity, assume $p$ is odd.  Then Brantner showed that $(E_n)_*\cLie(Z)$ is the free complete Hecke-Lie algebra on $(E_n)_*Z$:
$$ (E_n)_*\cLie(Z) = [\mc{H}^{\cLie} \otimes_{(E_n)_*} \cLie^{(E_n)_*}((E_n)_*Z)]^{\wedge}_I. $$
In the case of $n = 2$, the algebra $\mc{H}^{\cLie}$, and the Morava $E$-theory $(E_2)_*\Phi_2(S^{2i+1})$, was determined by Zhu \cite{Zhu}.

The first non-trivial differentials in the $K(n)$-local GSS are given by analogs of the formula (\ref{eq:rssdiff}).
In the case of $n = 1$ and $p$ odd, Kjaer used this to compute the $v_1$-periodic GSS in its entirety for $X$ a simply connected finite $H$-space \cite{Kjaer}.  By comparing his results with the work of Bousfield \cite{Bousfield99}, Kjaer established that for $p$ odd, all finite simply connected $H$-spaces are $\Phi_1$-good. This suggests that the higher chromatic analogs of the right-hand column of (\ref{eq:Quillen}) should be better behaved than the higher chromatic analogs of the left-hand column.  Progress on the study of $T(n)$-local Hopf algebras is being made in ongoing work of  Brantner, Hahn, Heuts, and Yuan, who have proposed that it may be the case that \emph{all} simply connected loop spaces are $\Phi_n$-good.

\bibliographystyle{amsalpha}
\bibliography{lie}

\end{document}